\documentclass[11pt,letterpaper]{article}

\usepackage{hyperref}

\usepackage{amssymb,amsmath}

\usepackage{epsf}
\usepackage{graphicx,amssymb,amsmath,color,latexsym} 

\usepackage{xcolor,wasysym,floatflt}

\newtheorem{thm}{Theorem}
\newtheorem{lem}{Lemma}

\newtheorem{prop}{Proposition}

\usepackage[none]{hyphenat} 

\makeatletter
\newcommand\footnoteref[1]{\protected@xdef\@thefnmark{\ref{#1}}\@footnotemark}
\makeatother
\def\qed{\hfill $\square$}

\begin{document}

\begin{center}
   \Large \textbf{Global asymptotic stability of nonconvex sweeping processes}  
 \vskip 0.3cm \large Lakmi Niwanthi Wadippuli, Ivan Gudoshnikov, Oleg Makarenkov
  \vskip.1cm \small Department of Mathematical Sciences, The University of Texas at Dallas, 800 West Campbell Road 
  Richardson, Texas 75080, USA
\end{center}

\noindent Abstract: Building upon the technique that we developed earlier for perturbed sweeping processes with convex moving constraints and monotone vector fields (Kamenskii et al, Nonlinear Anal. Hybrid Syst. 30, 2018), the present paper establishes global asymptotic stability of global and periodic solutions to perturbed sweeping processes with prox-regular moving constraint. Our conclusion can be formulated as follows:  closer the constraint to a convex one, weaker monotonicity is required to keep the sweeping process globally asymptotically stable. We explain why the proposed technique is not capable to prove global asymptotic stability of a periodic regime in a crowd motion model (Cao-Mordukhovich, DCDS-B  22, 2017). We introduce and analyze a toy model which clarifies the extent of applicability of our result. 


\section{Introduction}

Let  $t\mapsto C(t)$ be a set valued map which take nonempty closed values and  $f: \mathbb{R} \times \mathbb{R}^n \rightarrow \mathbb{R}^n$. Then the corresponding perturbed Moreau sweeping process is given as
\begin{equation}\label{spp}
-\dot{x} \in N\left(C(t),x\right)+f(t,x)
\end{equation}
where $N(C(t),\cdot)$ is normal cone to the set $C(t)$, given by
\begin{equation*}
N(C,x)=\{v\in\mathbb{R}^n : x\in {\rm proj}(x+\alpha v,C) \text{ for some } \alpha>0\}
\end{equation*}
and ${\rm proj}(x,C)$ is the set of points of $C$ closest to the point $x$.

\vskip0.2cm

\noindent We say an absolutely continuous function $x$ is a solution of sweeping process (\ref{spp}) on an interval $I\subset \mathbb{R}$ if $x(t)\in C(t)$ for each $t$ and $\dot{x}(t)$ satisfy (\ref{spp}) for a.e. $t\in I$.

\vskip0.2cm

\noindent Due to challenges from crowd motion modeling (Maury-Venel \cite{crowd}), the existence and uniqueness of a solution to nonconvex sweeping processes has being intensively studied.
The main problem of weakening the convexity of the set is the lack of continuity of the map $x \mapsto {\rm proj}(x,C)$ in general. Therefore, the concept of prox-regularity came to the study of sweeping processes. For the space $\mathbb{R}^n$, the set $C(t)$ is $\eta$-prox-regular, if $C(t)$ admits an external tangent ball with radius smaller than $\eta$ at each $x\in \partial C(t)$ (see Maury-Venel \cite[p.~150]{crowd}, Colombo and Monteiro Marques \cite[p.~48]{g-m}).

\vskip0.2cm

\sloppy

\noindent Colombo-Goncharov \cite{colombo}, Benabdellah \cite{benabdellah}, Colombo and Monteiro Marques \cite{g-m}, and Thibault \cite{thibault} studied the existence and uniqueness of solutions to non-perturbed sweeping processes with nonconvex prox-regular sets. Existence and uniqueness for  perturbed sweeping processes is considered in Edmond-Thibault \cite{edm1}, \cite{edm2}. A sweeping process with  prox-regular set values appeared in the context of crowd motion modeling in Maury-Venel \cite{crowd} along with numerical simulations. Cao-Mordukhovich \cite{corridor} illustrate their result for nonconvex sweeping process using crowd motion model of traffic flow in a corridor. Edmond-Thibault \cite{edm2}, Cao-Mordukhovich \cite{cao-crowd} studied optimal control problems related to a nonconvex perturbed sweeping process. Optimal control problem of convex sweeping process which is coupled with a differential equation was studied in Adam-Outrata \cite{adam} and the possibility of weakening the convexity to prox-regularity is mentioned there.

\vskip0.2cm

\noindent The problem of the existence of periodic solutions in sweeping processes with convex constraint was of interest lately, see e.g. Castaing and Monteiro Marques \cite[Theorem~5.3]{castaing}, Kunze \cite{kunze} and Kamenskii-Makarenkov \cite{kamenskii} and references therein.
\vskip0.2cm

\noindent In this paper we investigate stability of both arbitrary global solution and a periodic solution of sweeping processes (\ref{spp}) with prox-regular set-valued function $C(t)$. The existence of globally exponentially stable global and periodic solutions to (\ref{spp}) when $C(t)$ is convex-valued has been recently established in Kamenskii et al \cite{glob}. The central setting of \cite{glob} is  strong monotonicity of  $f$ in the sense that 
\begin{equation}\label{monot}
\langle f(t,x_1)-f(t,x_2),x_1-x_2 \rangle  \geq \alpha \lVert x_1-x_2 \lVert ^2, \quad \mbox{for all}\ t\in\mathbb{R},\ x_1,x_2\in\mathbb{R}^n,
\end{equation}
for some fixed $\alpha >0$.
A similar framework has been earlier used by  Heemels-Brogliato \cite{brogliato}, Brogliato \cite{brogliatoF} and Leine-van de Wouw \cite{leine} to prove incremental stability of sweeping process (\ref{spp}) with time-independent convex constraint.  The present paper, for the first time ever, takes advantage of property (\ref{monot}) in the context of prox-regular non-convex sets $C(t)$. 

\vskip0.2cm

\noindent The paper is organized as follows. The next section is devoted to the proof of the main result (Theorem~\ref{periodic}), which gives conditions for global asymptotic stability of a periodic solution to (\ref{spp}). The structure of our proof is motivated by the method of our paper \cite{glob}. Indeed, the existence of a global solution to (\ref{spp}) follows the lines of the proof of Theorem 2.1 in \cite{glob} since the proof is independent of the convexity of the set (the proof of Theorem~\ref{ex} is still given in Appendix for completeness). At the same time, additional assumptions, compared to \cite{glob} are still required. First of all, in order to use the hypomonotonicity of the prox normal cone,  we need $f(\cdot,x)$ to be globally bounded for each $ x\in \bigcup\limits_{t\in \mathbb{R}}C(t)$, 
additionally to the assumptions of Theorem 2.2 in \cite{glob}. Furthermore, to obtain contraction of solutions to sweeping process (\ref{spp}), a lower bound of constant $\alpha$ in (\ref{monot}) depending on prox-regularity constant of the set $C(t)$ is required (Theorem~\ref{uni}).  

\vskip0.2cm

\noindent Section 3 is devoted to examples that illustrate the main result. Though global stability of the sweeping process of crowd motion model of  Maury-Venel \cite{crowd} has been the main driving force behind this paper, it still remains an open question as we discuss in the Appendix.


\section{The main result}
Let $C:\mathbb{R} \to \mathbb{R}^n$ be a nonempty closed $\eta$-prox-regular set-valued function with  Lipschitz continuity 
\begin{equation}\label{LipC}
d_H(C(t_1),C(t_2))\le L_C|t_1-t_2|,\quad \text{for all}\ t_1,t_2\in\mathbb{R}, \text{ and for some }L_C>0,
\end{equation} 
where $d_H(C_1,C_2)$ is the Hausdorff distance between two closed sets $C_1,C_2 \subset \mathbb{R}^n$ given by
\begin{equation}\label{dH}
d_H(C_1,C_2)=\max\left\{  
\sup_{x\in C_2} {\rm dist}(x,C_1),\sup_{x\in C_1} {\rm dist}(x,C_2)
\right\}
\end{equation}
with ${\rm dist}(x,C)=\inf\left\{
|x-c|:c\in C
\right\}.$\\\\
And let $f: \mathbb{R}\times \mathbb{R}^n \to \mathbb{R}^n$ be such that for some $L_f>0$ 
\begin{equation}\label{Lipf}
\|f(t_1,x_1)-f(t_2,x_2)\|\le L_f\|t_1-t_2\|+L_f\|x_1-x_2\|,
\end{equation} 
$\text{for all }t_1,t_2\in\mathbb{R},\ x_1,x_2\in\mathbb{R}^n.$

\vskip0.2cm

\noindent Here we will be using the hypomonotonicity of the normal cone for $\eta$-prox-regular sets Edmond-Thibault \cite[p. 350]{edm2} which is given as
\begin{equation}\label{hypo}
\langle v-v' , x-x'\rangle \geq - \lVert x-x'\rVert ^ 2 
\end{equation}
$\text{for } v\in N(C,x), v'\in N(C,x') \text{ such that } \lVert v \rVert, \lVert v' \rVert \leq \eta.$

\vskip0.2cm

\noindent We will be using the following version of Gronwall-Bellman lemma Trubnikov-Perov \cite[Lemma~1.1.1.5]{trub} (see also Kamenskii et al \cite[lemma 6.1]{glob}) in our proofs. 

\begin{lem}\label{trub} {\bf (Gronwall-Bellman)} Let an absolutely continuous function $a:[0,T]\to\mathbb{R}$ satisfy
	\begin{equation*}\label{1.1.1.21}
	\dot{a}(t) \le \lambda a(t)+b(t),\qquad{\rm for\ a.e.\ }t\in[0,T],
	\end{equation*}
	where $b:[0,T]\to\mathbb{R}$ is an integrable function and $\lambda\in\mathbb{R}$ is a constant. Then
	$$
	a(t)\le e^{\lambda t}a(0)+\int\limits_0^t e^{\lambda (t-s)}b(s)ds,\qquad{\rm for\ all\ }t\in[0,T].
	$$
\end{lem}

\begin{thm}\label{ex}
	Let $ C: \mathbb{R} \to \mathbb{R}^n$ be a Lipschitz continuous function with constant $L_C$ and let $C(t)$ be nonempty, closed and $\eta$-prox-regular for each $t \in \mathbb{R}$. Let $f: \mathbb{R} \times \mathbb{R}^n \to \mathbb{R}^n$ satisfy Lipschitz condition (\ref{Lipf}).
	Then the sweeping process (\ref{spp}) has at least one solution defined on the entire $\mathbb{R}$. 
\end{thm}
\noindent The proof follows same steps as in the proof of Theorem~2.1 in \cite{glob}. But we include the proof in the Appendix for completeness of the paper.

\begin{thm}\label{uni}
	Let the conditions of Theorem~\ref{ex} hold and $L_C\geq 0$ is as given by Theorem~\ref{ex}.
	Let 
	\begin{equation}\label{unifbound}	
	\lVert f(t,x)\rVert \leq M_f, \mbox{ for all } t\in \mathbb{R},\  x\in \bigcup\limits_{t\in \mathbb{R}}C(t),
	\end{equation} 
	where $M_f\geq 0$ is a fixed constant. Assume (\ref{monot}) holds with
	\begin{equation}\label{inequal}
	\alpha >\dfrac{L_C+M_f}{\eta}.
	\end{equation}
	Then the sweeping process (\ref{spp}) has a unique solution $x_0$, defined on $\mathbb{R}$. Furthermore the global solution $x_0$ is globally exponentially stable. 
\end{thm}
\noindent {\bf Proof.} 
	We note that by Edmond-Thibault \cite[Proposition~1]{edm2} for a solution $x$ of (\ref{spp}) with the initial condition $x(\tau)=x_0$,
	\begin{equation*}
	\lVert \dot{x}(t) + f(t,x(t)\rVert \leq \lVert f(t,x(t)\rVert +L_C, \quad \text{ for } t > \tau.
	\end{equation*}
	Then with uniform boundedness of $f$ we have
	\begin{equation}\label{bdd}
	\lVert \dot{x}(t) + f(t,x(t)\rVert \leq M_f +L_C, \quad \text{ for } t > \tau.
	\end{equation}
	Now let $x_1, x_2$ be two solutions of (\ref{spp}) with initial conditions $x_1(\tau), x_2(\tau) \in C(\tau)$. Let $t\geq \tau$ such that $\dot{x}_1(t), \dot{x}_2(t)$ exist.
	\\	
	Since $$ -\dot{x}_1(t) - f(t,x_1(t))\in N(C(t),(x_1(t)), \quad -\dot{x}_2(t) - f(t,x_2(t)) \in N(C(t),(x_2(t)),$$
	by hypomonotonicity condition (\ref{hypo}) of the normal cone and by (\ref{bdd}) we have
	\begin{align*}
	\left\langle \frac{-\eta}{M_f+L_C}(\dot{x}_1(t)+f(t,x_1(t))) -\frac{-\eta}{M_f+L_C}(\dot{x}_2(t)+ f(t,x_2(t))) ,  x_1(t)-x_2(t)\right\rangle \\ 
	\geq -\lVert x_1(t)-x_2(t) \lVert ^2 .
	\end{align*}
	Then
	\begin{align*}
	\lVert x_1(t)-x_2(t) \lVert ^2  - \frac{\eta}{M_f+L_C}&\langle f(t,x_1(t)) -  f(t,x_2(t)) , x_1(t)-x_2(t)\rangle  \\ 
	& \geq \frac{\eta}{M_f+L_C}\langle \dot{x}_1(t) - \dot{x_2}(t) , x_1(t)-x_2(t)\rangle,
	\end{align*}
	and by (\ref{monot}),
	\begin{equation*}
	\frac{\eta}{M_f+L_C}\langle \dot{x}_1(t) -\dot{x_2}(t) , x_1(t)-x_2(t)\rangle \leq \lVert x_1(t)-x_2(t) \lVert ^2 - \frac{\eta\alpha}{M_f+L_C}\lVert x_1(t)-x_2(t) \lVert ^2.
	\end{equation*}
	Thus we have
	\begin{equation*}
	\langle \dot{x}_1(t) -\dot{x_2}(t) , x_1(t) - x_2(t)\rangle \leq \left(\frac{M_f+L_C}{\eta} -\alpha\right)\lVert x_1(t)-x_2(t) \lVert ^2,
	\end{equation*}
	i.e. 	\begin{equation*}
	\frac{d}{dt} \lVert x_1(t)-x_2(t) \lVert ^2 \leq \left(\dfrac{2(M_f+L_C)}{\eta} -2\alpha\right)\lVert x_1(t)-x_2(t) \lVert ^2. 
	\end{equation*}
	Let $\bar{\alpha} = \frac{1}{\eta}\left(M_f+L_C -\eta\alpha \right) $.
	Then by Gronwall-Bellman lemma (\ref{trub}), for $t>\tau$,
	\begin{equation*}
	\lVert x_1(t)-x_2(t) \lVert ^2 \leq e^{2\bar{\alpha}(t-\tau)} \lVert x_1(\tau)-x_2(\tau) \lVert^2,
	\end{equation*} 
	and so \begin{equation}\label{stab}
	\lVert x_1(t)-x_2(t) \lVert  \leq e^{\bar{\alpha}(t-\tau)} \lVert x_1(\tau)-x_2(\tau)\lVert, \quad \text{for }\ t>\tau. 
	\end{equation} 
	\\
	
	\noindent Let $x(t)$ be a global solution of (\ref{spp}) which exists by Theorem~\ref{ex}. Then (\ref{inequal}) guarantees that $\bar \alpha <0$  and that $x(t)$ is exponentially stable. It remains to observe that $x(t)$ is the only global solution. Indeed, let $\bar x(t)$ be another global solution. Then, for each $t \in \mathbb{R}$ we can pass to the limit as $\tau \to \infty$ in  (\ref{stab}), obtaining $\|x(t)-\bar x(t)\|\leq 0$, so $x=\bar x$. 
\qed\\\\
\noindent Now we give a theorem about periodicity of the unique global solution established in Theorem \ref{uni}. The proof follows the lines of Castaing and Monteiro Marques \cite[Theorem~5.3]{castaing}, but we include such a proof for completeness. 
\begin{thm}\label{periodic}
	The unique global solution $x_0$ which comes from Theorem~\ref{uni} is T-periodic, if both maps $t\mapsto C(t)$ and $t\mapsto f(t,x)$ are T-periodic. 
\end{thm}
\noindent {\bf Proof.}
	Note that $a\mapsto x_a(T)$ is a contraction mapping from $C(0)$ to $C(T)=C(0)$, where $x_a$ is the solution of (\ref{spp}) on $[0, T]$ with initial condition $x_a(0)=a \in C(0)$.  
	Indeed, by (\ref{stab}), for $a,b\in C(0)$,
	\begin{equation*}
	\lVert x_a(T)-x_b(T) \lVert  \leq e^{\bar{\alpha}T} \lVert a-b\lVert 
	\end{equation*} 
	where $\bar{\alpha}<0$. 
	
	\vskip0.2cm
	\noindent Then, since $a\mapsto x_a(T)$ is continuous on $C(0)$ (see Edmond-Thibault \cite[Proposition~2]{edm2}), by the contraction mapping principle on $C(0)$ (see Rudin \cite[p.220]{rudin}), there exists $\bar{x} : [0, T] \to C(0)$ such that $\bar{x}(0)=\bar{x}(T)$ and satisfies (\ref{spp}) on $[0, T]$. Since both $t\mapsto C(t)$ and $t\mapsto f(t,x)$ are $T$-periodic, we can extend $\bar{x}$ to a $T$-periodic solution defined on $\mathbb{R}$ by $T$-periodicity.
	
	\vskip0.2cm
	
	\noindent	Since the global solution $x_0$ given by Theorem~\ref{uni} is unique, we have the result.
\qed

\section{Example}
Let the vector field $f: \mathbb{R}\times\mathbb{R}^2  \to \mathbb{R}^2$ be given by 
\begin{equation}\label{fexample}
f(t,x):=\alpha x, \quad t\in \mathbb{R},\ x\in\mathbb{R}^2,
\end{equation}
where $\alpha>0$ is a fixed constant.
We define the moving set $C(t)$ using a function $b\in C^1(\mathbb{R},  \mathbb{R})$ which is bounded below by $\beta\geq 1$ and admits a global Lipschitz  constant $L_b$, i.e.
\begin{equation}\label{conb}
|b(t_1)-b(t_2)|\leq L_b|t_1-t_2|, \quad \text{for all } \ t_1,t_2\in \mathbb{R}.
\end{equation}
Define 
\begin{equation}
C(t):=\bar{B}_{1}  \bigcap  S(t), \quad S(t)=\left\lbrace x\in \mathbb{R}^2 : x^2_1 + \dfrac{x^2_2}{b(t)^2}\geq 1\right\rbrace.
\label{eq:C}
\end{equation}
where $\bar{B}_{1}$ is the closed ball of radius 1 and centered at $(-1.5,0)$.  
\begin{figure}[h]
	\begin{center}
		\includegraphics[width=0.9\textwidth]{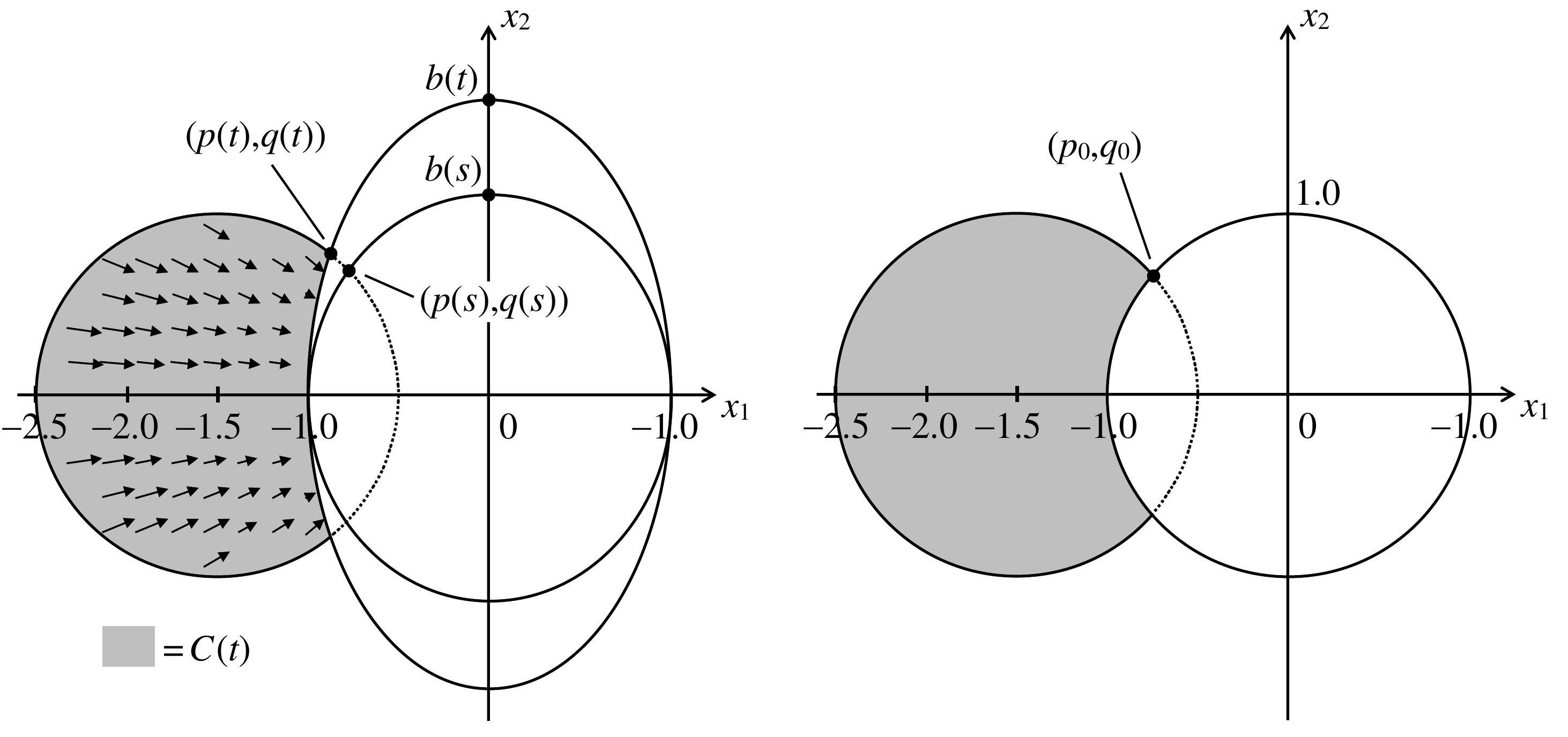}
	\end{center}
	\caption{Illustrations of the notations of the example. The closed ball centered at $(-1.5,0)$ is $\bar B_1$ and the white ellipses are the graphs of $S(t)$ for different values of the argument. The arrows is  the vector field of $\dot{x}=-\alpha x$.}\label{olegfig}
\end{figure}

\noindent In order to apply Theorem~\ref{uni}, we will now analyze: {\it i)} strong monotonicity and uniform boundedness of $f(t,x)$, {\it ii)} Lipschitz continuity of $C(t)$, {\it iii)} prox-regularity of $C(t)$.

\vskip0.2cm

\noindent {\it i) The monotonicity and boundedness of $f(t,x).$}  Since $\langle f(t,x)-f(t,y),x-y\rangle=\langle\alpha x-\alpha y,x-y \rangle=\alpha \|x-y\|^2$, $f$ is strongly monotone with constant $\alpha$ and bounded on $\bar{B}_1 \supset C(t)$ by $M_f=2.5\alpha$.\\

\noindent {\it ii) Lipschitz continuity of $C(t)$.} The boundary $\partial \bar B_1$ of $\bar B_1$ intersects the boundary $\partial S(t)$ of $S(t)$ at a unique point $(p(t),q(t))$ with  $q(t)\ge 0$. Since
$$
d_H(C(t),C(s))\le\|(p(t),q(t))-(p(s),q(s))\|
$$
(see Fig.~\ref{olegfig}),
we now aim at computing the Lipschitz constants of functions $p$ and $q.$ Since $b\in C^1(\mathbb{R},[1,\infty)),$ the implicit function theorem (see e.g. Zorich \cite[Sec.~8.5.4, Theorem~1]{zorich}) ensures that $p$ and $q$ are differentiable on $\mathbb{R}.$ Therefore, by the mean-value theorem (see e.g. Rudin \cite[Theorem~5.10]{rudin}),
\begin{equation}\label{oleg1}
d_H(C(t),C(s))\le\|(p'(t_p),q'(t_q))\|\cdot |t-s|,
\end{equation}
where $t_p,t_q$ are located between $t$ and $s.$  To compute $(p'(t_p),q'(t_q))$, we use the formula for the derivative of the implicit function (Zorich \cite[Sec. 8.5.4, Theorem~1]{zorich})
$$
(p'(t),q'(t))^T=-\left(F'_{(p,q)}\right)^{-1}(p(t),q(t),t)F'_t(p(t),q(t),t),
$$
applied with
$$
F(p,q,t)=\left(\begin{array}{c} 
(p+1.5)^2+q^2-1\\ 
p^2+\dfrac{q^2}{b(t)^2} -1
\end{array}\right).
$$
Since
$$
F'_{(p,q)}(p,q,t)=2\left(\begin{array}{cc} p+1.5 & q \\ p & \dfrac{q}{b(t)^2}\end{array}\right),\quad F'_t(p,q,t)=\left(\begin{array}{c} 0 \\ -2b(t)^{-3}b'(t)q^2\end{array}\right),
$$
we get the following formula for the derivatives $p'$ and $q'$
$$
\left(\begin{array}{c}  p'(t) \\ q'(t) \end{array}\right)=-\dfrac{1}{\dfrac{1}{b(t)^2}(p(t)+1.5)q(t)-p(t)q(t)}\left(\begin{array}{c} q(t) \\ -(p(t)+1.5)\end{array}\right)\dfrac{1}{b(t)^3}q(t)^2b'(t).
$$
Noticing that the properties $1+p(t)>0$ and $-p(t)b(t)^2>0$ imply 
$$
\dfrac{1}{b(t)\cdot(p(t)+1.5-p(t)b(t)^2)}\le 
\dfrac{1}{\beta\cdot(-p(t) b(t)^2)}\le\dfrac{1}{\beta^3 |p_0|},
$$
we conclude 
$$
|p'(t)|\le \dfrac{L_b}{\beta^3 |p_0|},
\qquad   |q'(t)|\le\dfrac{\color{black}L_b}{\beta^3 |p_0|},
$$
where $p_0$ is such that $p(t)\le p_0$ for all $t\in\mathbb{R}.$ Since $b(t)\ge 1,$ we can take $p_0$ as the abscissa of the intersection of $\partial \bar B_1$ with a unit circle centered at $0$, i.e.  
$$
p_0=-0.75, 
$$
see  Fig.~\ref{olegfig}.
Substituting these achievements to (\ref{oleg1}), we conclude
$$
d_H(C(t),C(s))\le \dfrac{{\color{black}4}L_b}{3\beta^3}|t-s|, 
$$
which gives $L_C=\dfrac{{\color{black}4}L_b}{3\beta^3}$ for the Lipschitz constant of $t\mapsto C(t)$.

\

\noindent {\it iii) The constant $\eta$ in $\eta$-prox-regularity of $C(t).$} 
We recall that $C(t)$ is $\eta$-prox-regular if $C(t)$ admits an external tangent ball with radius smaller than $\eta$ at each $x\in \partial C(t)$ (see Maury and Venel \cite{crowd}, Colombo and Monteiro Marques \cite{g-m}). The points of $\partial C(t)\backslash \partial S(t)$ admit an external tangent ball of any radius. Therefore, to find $\eta$, which determines $\eta$-prox-regularity of $C(t)$, it is sufficient to  focus on the points of $\partial C(t)\cap \partial S(t)$. That is why, for a fixed $t\in\mathbb{R},$  
we can choose $\eta$ as the minimum of the radius of curvature through 
$x\in\partial C(t)\cap \partial S(t)$, see 
e.g. Lockwood \cite[p.~193]{curvature}. 

\vskip0.2cm

\noindent Let us fix $t\in\mathbb{R}$ and use the parameterization $P(\phi)=(-\cos \phi,b(t)\sin \phi)$, $\phi\in \left[-\frac{\pi}{2},\frac{\pi}{2}\right],$ for the left-hand side of the ellipse $x^2 + \frac{y^2}{b(t)^2}= 1$. Then,  the radius of curvature $R(\phi)$ of  $\partial C(t)\cap \partial S(t)$ at $P(\phi)$  is (see Lockwood \cite[p.~xi,\ p.~21]{curvature}) 
\[
R(\phi)=\dfrac{1}{b(t)}(\sin^2 \phi+b(t)^2\cos^2\phi)^{\frac{3}{2}}=\dfrac{1}{b(t)}\left(b(t)^2+(1-b(t)^2)\sin^2(\phi)\right)^\frac{3}{2}.
\] 
Observe that $R$ decreases when $|\phi|$ increases from 0 to $\dfrac{\pi}{2}$. 

\begin{figure}[h]
	\begin{center}
		\includegraphics[width=0.45\textwidth]{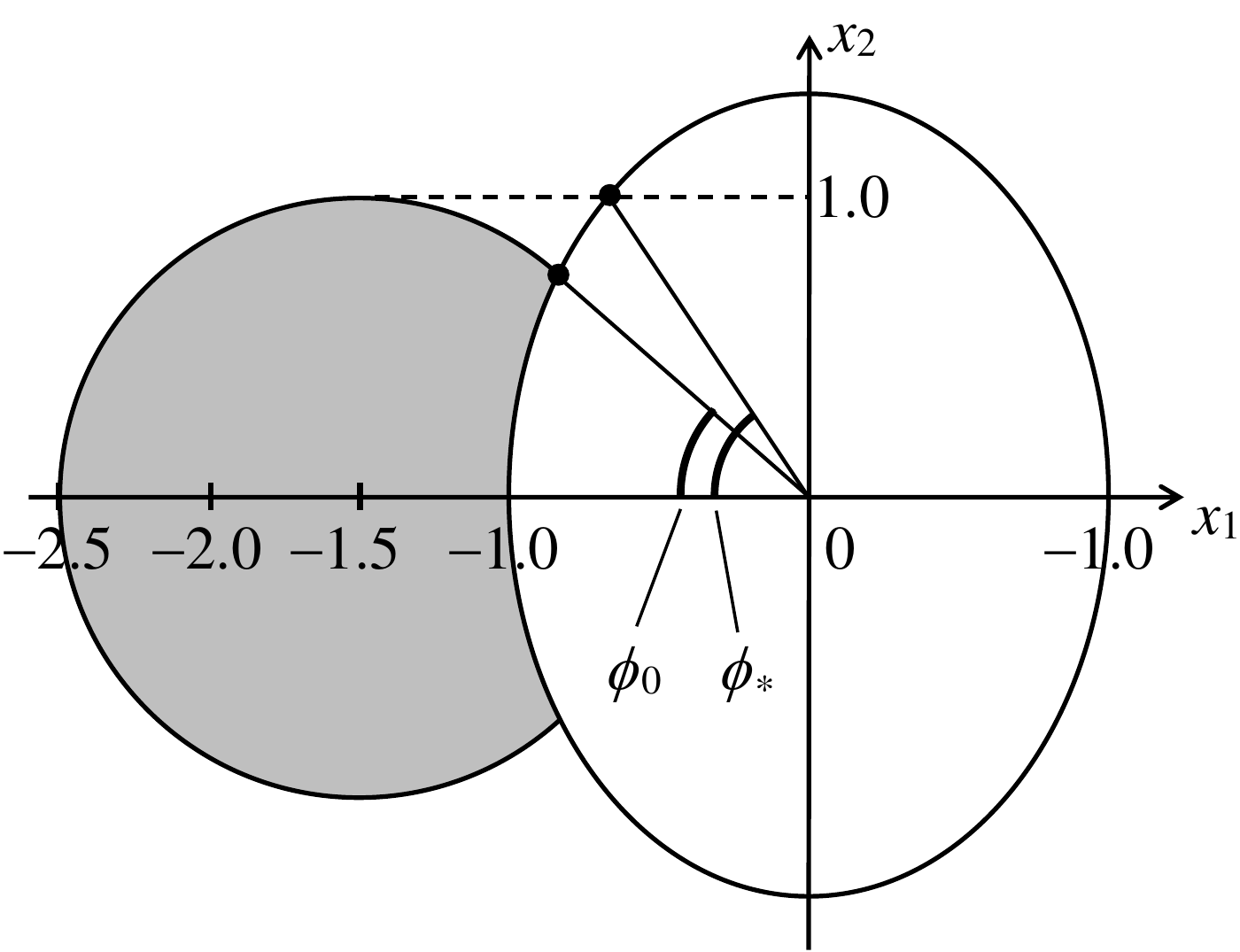}
	\end{center}
	\caption{The parameters $\phi_0$ and $\phi_*.$}\label{phifig}
\end{figure}

\noindent Therefore, the minimum curvature of $\partial C(t)\cap \partial S(t)$ is attained at the point $(p(t),q(t))$ as defined in ii). Let $\phi_0$ be such that $P(\phi_0)=(p(t),q(t))$ and let $\phi_*>0$ be such  that
the second component $P_2(\phi_*)$ of $P(\phi_*)$ equals 1, which exists because $b(t)\ge 1$ (see Fig.~\ref{phifig}). Since $q(t)\le 1$, we have $\phi_0\le \phi_*,$ and since $\phi\to R(\phi)$ decreases as $|\phi|$ increases, we have 
$$
R(\phi_0)\ge R(\phi_*).
$$	
Since $P_2(\phi_*)=1$ implies $b(t)\sin\phi_*=1$, we have $\sin\phi_*=\dfrac{1}{b(t)}$ and so
\begin{eqnarray*}
	R(\phi_0)&\ge& \dfrac{1}{b(t)}\left(\dfrac{1}{b(t)^2}+b(t)^2\left(1-\dfrac{1}{b(t)^2}\right)\right)^{\frac{3}{2}}=\dfrac{1}{b(t)}\cdot \dfrac{(1+b(t)^4-b(t)^2)^{\frac{3}{2}}}{b(t)^3}=\\
	&&=\left(b(t)^{-\frac{8}{3}}+b(t)^{\frac{4}{3}}-b(t)^{-\frac{2}{3}}\right)^\frac{3}{2}\ge  \left(b(t)^{\frac{4}{3}}-b(t)^{-\frac{2}{3}}\right)^\frac{3}{2}.
\end{eqnarray*}
Noticing that the function $b\mapsto \left(b^{\frac{4}{3}}-b^{-\frac{2}{3}}\right)^\frac{3}{2}$ increases on $[1,\infty),$ we finally conclude
$$
R(\phi_0)\ge \left(\beta^{\frac{4}{3}}-\beta^{-\frac{2}{3}}\right)^\frac{3}{2}.
$$
Therefore, $C(t)$ is $\eta$-prox-regular with $\eta=\left(\beta^{\frac{4}{3}}-\beta^{-\frac{2}{3}}\right)^\frac{3}{2}.$

\noindent Substituting the values of $M_f,$ $L_C,$ and $\eta$ into formula (\ref{inequal}), we get the following statement. 
\begin{prop}\label{propexample} Let $\alpha>0$ be an arbitrary constant and $b\in C^1(\mathbb{R},[\beta,\infty))$ with some $\beta\ge 1$ and Lipschitz condition (\ref{conb}). If
	$$
	\alpha>\dfrac{\dfrac{{\color{black} 4}L_b}{3\beta^3}+\dfrac{5}{2}\alpha}{\left(\beta^\frac{4}{3}-\beta^{-\frac{2}{3}}\right)^\frac{3}{2}},
	$$
	then, 
	the global solution
	\[
	x(t) = (-1,0), \quad  t\in \mathbb{R},
	\] of the sweeping process (\ref{spp}) with $C(t)$ and $f(t,x)$ given by (\ref{eq:C}) and (\ref{fexample}), 
	is globally asymptotically stable.
\end{prop}

\noindent As noticed earlier, $b\mapsto \left(b^{\frac{4}{3}}-b^{-\frac{2}{3}}\right)^\frac{3}{2}$ increases on $[1,\infty),$ so that the condition of Proposition~\ref{propexample} is a lower bound on $\beta.$

\section{Conclusion}
In this paper we  proved the existence of at least one global solution to a nonconvex sweeping process with Lipschitz right-hand-sides. The uniqueness and  exponential stability of the solution follows when the vector field of the sweeping process  is uniformly bounded, strongly monotone and the prox-regularity constant of the moving constraint is not too small. We further proved that the unique global solution is periodic when the right-hand-sides of the sweeping process are periodic in time. 

\vskip0.2cm

\noindent Following the lines of Kamenskii et al \cite{glob}, the ideas of the present work can be extended to almost periodic solutions and to sweeping processes with small non-monotone ingredients.

\vskip0.2cm

\noindent We show in Appendix that the estimate for the prox-regularity constant in  Maury-Venel \cite[Proposition~2.15, Proposition~2.17]{crowd} does not agree with inequality (\ref{inequal}), making our main result inapplicable to the model of \cite{crowd}. At the same time, we analyze a toy example where we document how  applicability or inapplicability of our result is linked to the parameters of sweeping process. 

\vskip0.2cm

\noindent The ultimate conclusion of the paper is as follows: closer the constraint to a convex one, weaker monotonicity is required to keep the sweeping process globally asymptotically stable.
\section{Appendix}

\subsection {\bf 
	Proof of Theorem~\ref{ex}.}  
Let $\{\xi_n\}_{n=1}^{\infty} \subset \mathbb{R}^n $ such that $\xi_n \in C(-n)$ for each $n\in \mathbb{N}$. Define \begin{equation*}
x_n(t) = \begin{cases}
x(t,-n,\xi_n) \quad \text{ if } t\geq-n  
\\ \xi_n \qquad \qquad\quad \text{ if } t < -n
\end{cases}
\end{equation*}
where  $x(t,-n,\xi_n)$ is the solution of (\ref{spp}) with initial condition $x(-n,-n,\xi_n)=\xi_n $ for $n\in \mathbb{N}$. By Edmond-Thibault \cite[Theorem~1]{edm2}, for each $n\in\mathbb{N}$, $x_n$ has the same Lipschitz constant $L_k>0$ on each interval $[-k,k]$ for each $k\in\mathbb{N}$. 
\\Let denote $x_n^0(t) = x_n(t)$ on $\mathbb{R}$ for each $n\in \mathbb{N}$. Then by Arzela-Ascoli theorem there exists a subsequence $\{x_n^{k}\}_{n=1}^{\infty} \subset \{x_n^{k-1}\}_{n=1}^{\infty}$ which converges uniformly on $[-k,k]$ for each $k\in\mathbb{N}$.
\\Now let define $\bar{x}_n=x^n_n$ on $\mathbb{R}$ for each $n\in \mathbb{N}$. Then $\{\bar{x}_n\}_{n\in \mathbb{N}}$ converges uniformly on $[-k, k]$ for each $k\in \mathbb{N}$. Let $x_0(t) := \lim\limits_{n\to \infty}\bar{x}_n(t)$.
\\\\Now let's show that $x_0$ is a solution of (\ref{spp}). 
\\Let $\bar{x}$ be a solution of (\ref{spp}) with initial condition $\bar{x}(\tau)=x_0(\tau)$. Assume $ \bar{x}(t_0) \neq x_0(t_0)$ for some $t_0 > \tau$. i.e. $ \lim\limits_{n\to\infty} \bar{x}_n(t_0) \neq \bar{x}(t_0)$. Then there exist $ \varepsilon_0 > 0 $ and  for each $n\in \mathbb{N}$,  $m_n>n$ such that $ \lVert \bar{x}_{m_n}(t_0) - \bar{x}(t_0)\lVert \geq \varepsilon_0 $. 
\\Then by continuously dependence of solution on the initial condition (see Edmond-Thibault \cite[Proposition~2]{edm2}), there exists $ \delta > 0$ such that if $ \lVert \bar{x}(\tau) - \bar{x}_n(\tau) \lVert < \delta$ then $ \lVert \bar{x}(t) - \bar{x}_n(t) \lVert < \varepsilon_0 $ for $n \in \mathbb{N}$ with $-n>\tau$ on $[\tau,t_0]$. 
\\But since $\bar{x}(\tau)=x_0(\tau) = \lim\limits_{n\to\infty} \bar{x}_n(\tau) $, there exists $N \in \mathbb{N}$ such that $ \lVert \bar{x}(\tau) - \bar{x}_n(\tau) \lVert < \delta $ for each $n>N$. Then $ \lVert \bar{x}(t) - \bar{x}_n(t) \lVert < \varepsilon_0 $ for $n >N$ on $[\tau,t_0]$. This contradicts $ \lim\limits_{n\to\infty} \bar{x}_n(t_0) \neq \bar{x}(t_0)$. Therefore $ \bar{x}(t) = x_0(t)$ for each $t\geq \tau$. Hence $x_0$ is a solution of (\ref{spp}).\\
The global boundedness of $x_0$ follows from the boundedness of $C$ on $\mathbb{R}$ and $x_0(t)\in C(t)$ for each $t\in \mathbb{R}$.
\qed

\subsection{The crowd motion model}
We give a brief introduction into the model by Maury-Venel \cite{crowd}, before we explain the inapplicability of  Theorem~\ref*{uni} in this model.
\\Consider $N$ people whose positions are given by $x=(x_1, x_2, \dots , x_N) \in \mathbb{R}^{2N}$, where each person is  identified as a disk with center $x_i \in \mathbb{R}^2$ and radius $r$.
\\ By avoiding overlapping of people, the set of feasible configurations is defined as
\begin{equation}\label{crowdC}
C=\{ x\in \mathbb{R}^{2N}:\lVert x_i-x_j \rVert -2r \geq 0 \text{ for all 
} i<j\}.
\end{equation}
Now let $U(x) = \left(U_1(x), U_2(x), \cdots U_N(x) \right)$ be the spontaneous velocity of each person at the position $x$ , i.e. $U_i(x)$ is the velocity that $i$-th person would have in the absence of other people. 
\\Since the aim of Maury-Venel \cite{crowd} is to have a model that describes people in a highly packed situation, the actual velocity of a person is defined to be closest to the spontaneous velocity. So the actual velocity is computed as the projection of the spontaneous velocity onto the set of feasible velocities. This gives the sweeping process
\begin{equation}\label{crowdSP}
\begin{cases}
-\dot{x} \in N(C,x) - U(x) 
\\x(0) = x_0 \in C.
\end{cases}
\end{equation}
\\Let's consider the situation where there are only two people. Then by Maury-Venel \cite[Proposition~2.15]{crowd}, the set $C$ in (\ref{crowdC}) is $\eta$-prox regular with $\eta = r\sqrt{2}$. Let's take $U(x)=-x$. Viewing (\ref{crowdSP}) as (\ref{spp}), we get $\alpha=1$ in (\ref{monot}). 
\\Then the condition (\ref{inequal}) of Theorem~\ref{uni} takes the form $\sqrt{2}r >L_C+M_f$, where $L_C=0$ (because $C$ in (\ref{crowdSP}) doesn't depend on $t$) and $M_f$ satisfies $\lVert f(t,x) \rVert = \lVert x \rVert \leq M_f$ for each $x\in C$. Therefore (\ref{inequal}) implies $M_f < \sqrt{2}r$.\\
On the other hand, since $\lVert (0,-r)-(0,r) \rVert =2r$, we have $(0,-r,0,r) \in C$ and so $M_f$ must verify $M_f \geq \lVert (0,-r,0,r) \rVert = \sqrt{2}r$. 
\\Therefore Theorem~\ref{uni} does not apply.

\end{document}